\documentclass[12pt]{article}
\usepackage{latexsym, amsfonts, amscd, amssymb, verbatim, amsmath,enumerate}

\begin{document}

\title{\textbf{On V-Semirings and Semirings all of whose Cyclic Semimodules
are Injective}}
\author{\textbf{\ }J. Y. Abuhlail$^{1}$, S. N. Il'in$^{2}$, Y. Katsov$^{3}$,
T.G. Nam$^{4}$ \\
{\footnotesize abuhlail@kfupm.edu.sa; Sergey.Ilyin@kpfu.ru;
katsov@hanover.edu; }\\
{\footnotesize trangiangnam05@yahoo.com;}\\
$^{1}${\footnotesize Department of Mathematics }$\&${\footnotesize \
Statistics, Box }$\#5046,${\footnotesize \ KFUPM, 31261 Dhahran, KSA }\\
$^{2}${\footnotesize Lobachevsky Institute of Mathematics and Mechanics}\\
{\footnotesize Kazan (Volga Region) Federal University, 18, Kremlevskaja
St., 420008 Kazan, }\\
{\footnotesize Tatarstan, Russia}\\
$^{3}${\footnotesize Department of Mathematics}\\
{\footnotesize Hanover College, Hanover, IN 47243--0890, USA}\\
$^{4}${\footnotesize Institute of Mathematics, VAST, 18 Hoang Quoc Viet, 10307,
Hanoi, Vietnam}\\
Dedicated to the memory of Professor L. A. Skornjakov}
\date{}
\maketitle

\begin{abstract}
In this paper, we introduce and study V- and
CI-semirings---semirings all of whose simple and cyclic,
respectively, semimodules are injective. We describe V-semirings for
some classes of semirings and establish some fundamental properties
of V-semirings. We show that all Jacobson-semisimple V-semirings are
V-rings. We also completely describe the bounded distributive
lattices, Gelfand, subtractive, semisimple, and anti-bounded,
semirings that are CI-semirings. Applying these results, we give
complete characterizations of congruence-simple subtractive and
congruence-simple anti-bounded CI-semirings which solve two earlier
open problems for these classes of CI-semirings.

\textbf{2010} \textbf{Mathematics Subject Classifications}: Primary 16Y60,
16D99, 06A12; Secondary 18A40, 18G05, 20M18

\end{abstract}

\footnotetext{%
The fourth author is partially supported by Vietnam National Foundation for
Science and Technology Development (NAFOSTED)}

\section{Introduction}

In the modern homological theory of modules over rings, the results
characterizing rings by properties of modules and/or suitable categories of
modules over them are of great importance and sustained interest (for a good
number of such results one may consult, for example, \cite{AF1997}, \cite%
{Lam1999}, and \cite{Wis1991}). Inspired by this, during the last three
decades quite a few results related to this genre have been obtained in
different nonadditive settings. Just to mention some of these settings, a
very valuable collection of numerous interesting results\ on
characterizations of monoids by properties and/or by categories of acts over
them, \textit{i.e.}, results in the so called \textit{homological
classification of monoids}, can be found in \cite{KKM2000}; and, for results
in the \textit{homological classification of distributive lattices}, another
non-additive setting, one may consult the survey \cite{Fof1982}.

Moreover, quite recently there was obtained a number of interesting and
important homological results in a more general, and gaining increasing
interest, non-additive setting --- results in the \textit{homological
classification/characterization of semirings} (see, for instance, \cite%
{Kat1997}, \cite{Kat2004}, \cite{Ili2006}, \cite{Ili2008},
\cite{Ili2010}, \cite{KN2011}, \cite{KNT2011}, \cite{Ili2012}, and
the papers cited in them). One may clearly notice a growing interest
in developing algebraic and homological theories of semirings and
semimodules motivated by their numerous connections with, and
applications in, different branches of mathematics, computer
science, quantum physics, and many other areas of science (see, for
example, \cite{Gla2002}). As algebraic objects, semirings are
certainly the most natural generalization of such, at first glance
different, algebraic systems as rings and bounded distributive
lattices. Thus, investigating semirings and their representations,
one should undoubtedly use methods and techniques of the ring,
lattice and semigroup theories, as well as diverse techniques and
methods of categorical and universal algebra. The wide variety of
the algebraic techniques involved in studying semirings, and their
representations/semimodules, perhaps explains why the research on
categorical and homological aspects of theory of semirings and
semimodules is still behind that for rings and monoids. In light of
this, presenting some new, important and \ interesting in our view,
nontrivial at all, homological considerations, results and
techniques to the problems of the homological
characterization/classification of semirings, as well as motivating
an interest to this direction of research, is a main goal of our
paper.

In this paper, we introduce and study semirings with two classes of
injective semimodules over them: \textit{V-semirings} \cite{Ili2012} ---
semirings all of whose simple semimodules are injective; and \textit{%
CI-semirings} --- semirings all of whose cyclic semimodules are injective.
The investigation and classification of such semirings serves as a
fundamental basis to obtain further structural insight of congruence-simple
semirings. The paper is organized as follows.

For the reader's convenience, all subsequently necessary and important
notions and facts on semirings and semimodules that cannot be found in \cite%
{Gol1999} and/or \cite{HW1996} are collected in Section 2.

In Section 3, together with constructing some new examples of
noncommutative V-semirings, we also characterize V-semirings within
important classes of semirings and establish some fundamental
properties of V-semirings. Among other results of this section, we
single out the following central ones: we give a complete
description of semisimple V-semirings (Theorem 3.12); it is shown
that the Jacobson-semisimple V-semirings are just the V-rings
(Theorem 3.14); it is established that for a semiring to be a
V-semiring is a Morita invariant property (Theorem 3.9).

In Section 4, among the main results of the paper are the following ones: we
describe the bounded distributive lattices, the Gelfand semirings \cite[p. 56%
]{Gol1999}, the left subtractive semirings, and the semisimple
semirings, that are CI-semirings (Theorem 4.3, Theorem 4.6, Theorem
4.7, and Theorem 4.10, respectively); we give a complete
characterization of anti-bounded CI-semirings (Theorem 4.20),
essentially generalizing B. Osofsky's celebrated characterization of
semisimple rings \cite{Oso1964} (see also \cite[Theorem
1.2.9]{Lam2001} and \cite[Corollary 6.47]{Lam1999}); and applying
Theorem 4.7 and Theorem 4.20, we give a complete description of
congruence-simple subtractive CI-semirings (Corollary 4.8) and
congruence-simple anti-bounded CI-semirings (Corollary 4.21),
respectively.

Finally, all notions and facts of categorical algebra, used here without any
comments, can be found in \cite{Mac1971}; for notions and facts from
semiring theory and universal algebra we refer to \cite{Gol1999} and \cite%
{Gra1979}, respectively.{\ {\ }}

\section{Preliminaries}

\noindent\ \ \ \ \ \ \textbf{2.1 }Recall \cite{Gol1999} that a \emph{%
semiring\/} is a datum $(S,+,\cdot ,0,1)$ such that the following conditions
are satisfied:\smallskip

(1) $(S,+,0)$ is a commutative monoid with identity element $0$;

(2) $(S,\cdot ,1)$ is a monoid with identity element $1$;

(3) Multiplication is distributive over addition from both sides;

(4) $0s=0=s0$ for all $s\in S$.\smallskip\

A semiring that is not a ring we call a \emph{proper semiring}.\smallskip\

As usual, a \emph{left\/} $S$-\emph{semimodule} over the semiring $S$ is a
commutative monoid $(M,+,0_{M})$ together with a scalar multiplication $%
(s,m)\mapsto sm$ from $S\times M$ to $M$ which satisfies the following
identities for all $s,s^{^{\prime }}\in S$ and $m,m^{^{\prime }}\in M$%
:\medskip

(1) $(ss^{^{\prime }})m=s(s^{^{\prime }}m)$;

(2) $s(m+m^{^{\prime }})=sm+sm^{^{\prime }}$;

(3) $(s+s^{^{\prime }})m=sm+s^{^{\prime }}m$;

(4) $1m=m$;

(5) $s0_{M}=0_{M}=0m$.\medskip

\emph{Right semimodules\/} over $S$ and homomorphisms between semimodules
are defined in the standard manner. And, from now on, let $\mathcal{M}$ be
the variety of commutative monoids, and $\mathcal{M}_{S}$ and $_{S}\mathcal{M%
}$ denote the categories of right and left $S$-semimodules, respectively,
over a semiring $S$.\smallskip\

\textbf{2.2 }An element $\infty \in M$ of an $S$-semimodule $M$ is \emph{%
infinite} if $\infty +m=\infty $ for every $m\in M$; and $K\leq _{\text{ }%
S}M $ means that $K$ is an $S$-subsemimodule of $M$. Also, we will use the
following subsets of the elements of an $S$-semimodule $M$ :%
\begin{equation*}
\begin{tabular}{lll}
$I^{+}(M)$ & $:=$ & $\{m\in M\,|\,m+m=m\};$ \\
$Z(M)$ & $:=$ & $\{z\in M\mid z+m=m\text{ for some }m\in M\};$ \\
$V(M)$ & $:=$ & $\{m\in M\mid m+m^{\prime }=0\text{ for some }m^{\prime }\in
M\};$ \\
$A(M)$ & $:=$ & $\{m\in M\mid $ $m+m^{\prime }+m=m$ for some $m^{\prime }\in
M\}.$%
\end{tabular}%
\end{equation*}%
\smallskip

For a semimodule $M$ $\in |_{S}\mathcal{M}|$, it is obvious that $%
I^{+}(M)\cap V(M)=\{0\}$, and $I^{+}(M)\leq _{\text{ }S}Z(M)\leq _{\text{ }%
S}M$. Moreover, if $S$ is an additively regular semiring, then it is easy to
see that $I^{+}(M)$ is a left $I^{+}(S)$-semimodule, as well as the
subsemimodule $V(M)$ is the largest $S$-module in $M$.

A left $S$-semimodule $M$ is \emph{zeroic} (\emph{zerosumfree,} \emph{%
additively idempotent, additively regular}) if $Z(M)=M$ ($V(M)=0,$ $%
I^{+}(M)=M,$ $A(M)=M$). In particular, a semiring $S$ is \emph{zeroic} (%
\emph{zerosumfree}, \emph{additively idempotent}, \emph{additively regular})
if $_{S}S$ $\in |_{S}\mathcal{M}|$ is a zeroic (zerosumfree, additively
idempotent, additively regular) semimodule. For example, the \emph{Boolean
semiring} $\mathbf{B}=\{0,1\}$ is a commutative, zeroic, zerosumfree,
additively idempotent semiring in which $\infty =1$.\smallskip

\textbf{2.3} A subsemimodule $K\leq _{S}M$ of a semimodule $M$ is (\emph{%
strong}) \emph{subtractive }if ($m+m^{\prime }\in K\Rightarrow m,m^{\prime
}\in K$) $m,m+m^{\prime }\in K\Rightarrow m^{\prime }\in K$ for all $%
m,m^{\prime }\in M$. A left $S$-semimodule $M$ is \emph{subtractive} if it
has only subtractive subsemimodules. A semiring $S$ is \emph{left subtractive%
} if $S$ is a subtractive left semimodule over itself. \emph{\ }

As usual (see, for example, \cite[Chapter 17]{Gol1999}), if $S$ is a
semiring, then in the category $_{S}\mathcal{M}$, a \emph{free}
(left) semimodule $\sum_{i\in I}S_{i},S_{i}\cong $ $_{S}S$, $i\in
I$, with a basis set $I$ is a direct sum (a coproduct) of $|I|$
copies of $_{S}S$; a semimodule $P\in |_{S}\mathcal{M}|$ is
\textit{projective} if it is a
retract of a free semimodule; a semimodule $F\in |_{S}\mathcal{M}|$ is \emph{%
flat} iff the functor $-\otimes _{S}F:\mathcal{M}_{S}\longrightarrow
\mathcal{M}$ preserves finite limits iff $F$ is a filtered (directed)
colimit of finitely generated free (projective) semimodules \cite[Theorem
2.10]{Kat1997}; a semimodule $M\in |_{S}\mathcal{M}|$ is \emph{finitely
generated}\textit{\ (\emph{cyclic})} iff $M$ is a homomorphic image of a
free left $S$-semimodule with a finite basis (a homomorphic image of $_{S}S$%
); a semimodule $M\in |_{S}\mathcal{M}|$ is \emph{injective} if for any
monomorphism $\mu :A\rightarrowtail B$ of left $S$-semimodules $A$ and $B\ $%
and every homomorphism $f\in $ $_{S}\mathcal{M(}A,M)$, there exists a
homomorphism $\widetilde{f}\in $ $_{S}\mathcal{M(}B,M)$ such that $%
\widetilde{f}\mu =f$.\smallskip

\textbf{2.4} \emph{Congruences }on an $S$-semimodule $M$ are defined in the
standard manner, and $\mathrm{Cong}(M)$ denotes the set of all congruences
on $M$. This set is non-empty since it always contains at least two
congruences---the \emph{diagonal congruence} $\vartriangle _{M}:=$ $\{(m,m)$
$|$ $m\in M$ $\}$ and the \emph{universal congruence }$M^{2}:=\{(m,n)$ $|$ $%
m,n\in M$ $\}$. Any \ subsemimodule $L\leq _{\text{ }S}M$ of an $S$%
-semimodule $M$ induces a congruence $\equiv _{L}$ on $M$, known as the
\emph{Bourne congruence}, by setting $m\equiv _{L}m^{\prime }$ iff $%
m+l=m^{\prime }+l^{\prime }$ for some $l,l^{\prime }\in L$; and $M/L$
denotes the factor $S$-semimodule $M/\equiv _{L}$ having the canonical $S$%
-surjection $\pi _{L}:M\longrightarrow M/L$.

Furthermore, a nonzero $S$-semimodule $M$ is \emph{simple
}(\emph{atom, s-simple}) if\linebreak
%
$\mathrm{Cong}(M)=\{\vartriangle _{R},M^{2}\}$ ($M$ has no
nonzero proper $S$-subsemimodules, $M$ has no nonzero proper subtractive $S$%
-subsemimodules). The following observations will prove to be
useful.\smallskip

\noindent \textbf{Proposition 2.5}\label{simple} \textit{For a nonzero }$S$%
\textit{-semimodule }$M$\textit{\ the following statements are true:}

\textit{\noindent (1) }$M$\textit{\ is simple iff every nonzero semimodule}$%
\ $\textit{\ homomorphism }$f:M\longrightarrow N$\textit{\ is injective;}

\textit{\noindent (2) If }$M$\textit{\ is simple, then }$M$\textit{\ is
s-simple.\smallskip }

\noindent \textit{Proof }(1) $\Longrightarrow $. Let $f:M\longrightarrow N$
be a nonzero morphism and $\equiv _{f}$ the congruence on $M$ defined by $%
m\equiv _{f}m^{\prime }$ iff $f(m)=f(m^{\prime })$. Since $M$ is a simple
semimodule and $f$ is a nonzero homomorphism, it is easy to see that $\equiv
_{f}\ =$ $\vartriangle _{M}$ and, hence, $f$ is injective.

$\Longleftarrow $. It is obvious.

(2) See \cite[Proposition 15.13]{Gol1999}.\textit{\ \ \ \ \ \ }$_{\square }$

\section{On V-semirings}

\qquad Generalizing the well known for rings notions and following \cite%
{Ili2012}, we call a semiring $S$ a \emph{left }(\emph{right})\emph{%
V-semiring} if every simple left (right) $S$-semimodule is injective; and an
$S$-semimodule $M$ is called an \emph{essential extension} of an

\noindent $S$-subsemimodule $L\leq _{\text{ }S}M$, $i:$ $L\rightarrowtail M$%
, if for every semimodule homomorphism $\gamma :M\longrightarrow N$, the
homomorphisms $\gamma i$ and $\gamma $ are simultaneously injective. The
following characterization of V-semirings will prove to be useful.\medskip\

\noindent \textbf{Theorem 3.1} \label{char-V}\emph{(\cite[Theorem 2.10]%
{Ili2012})} \textit{The following statements for a semiring }$S$\textit{\
are equivalent:}

\textit{(1) }$S$\textit{\ is a left (right)V-semiring;}

\textit{(2) Every essential extension of each simple left (right) }$S$%
\textit{-semimodule }$M$\textit{\ coincides with }$M$\textit{;}

\textit{(3) }$S\cong R\oplus Z$\textit{, where }$R$\textit{\ is a left
(right) V-ring and }$Z$\textit{\ is a zeroic left (right) V-semiring;}

\textit{(4) Every quotient semiring of }$S$\textit{\ is a left (right)
V-semiring.\medskip }

In particular, one may easily see that item (2) of this characterization
implies\medskip

\noindent \textbf{Corollary 3.2 } \textit{Finite direct products of left
(right) V-semirings are left (right) V-semirings.}\textbf{\medskip }

In this section, in addition to constructing some new examples of
noncommutative V-semirings, we also characterize V-semirings within
important classes of semirings, and establish some fundamental properties of
V-semirings.

From Theorem 3.1 it follows that for a zerosumfree V-semiring to be zeroic
is, in general, only a necessary condition. However, it is also a sufficient
condition if the semiring has only two trivial strong left (right)
ideals.\medskip

\noindent \textbf{Proposition 3.3}\label{zeroic} A \textit{zerosumfree
semiring }$S$ \textit{possessing only two strong left (right) ideals is a
left (right) V-semiring iff it} \textit{is a zeroic semiring.\smallskip }

\noindent \textit{Proof} $\Longrightarrow $. It follows immediately from
Theorem 3.1 (3).

$\Longleftarrow $. Let $M$ be an arbitrary simple left $S$-semimodule. By
\cite[Proposition 1.2]{Ili2012}, $(M,+,0)$ is either a group or an
idempotent monoid. Consider these two cases.

Let\ $(M,+,0)$ be a group. Since $S$ is zeroic, there exists $z\in S$ such
that $1+z=z$, and, hence, $zm=(1+z)m=1m+zm$ and $m=1m=0_{M}$ for every $m\in
M$, what contradicts $M\neq 0$.

Now, let $(M,+,0)$ be an idempotent monoid. For every $0\neq m\in M$,$\ $the
left annihilator $(0:_{S}m):=\{$ $s\in S$ $|$ $sm=0_{M}$ $\}$ is obviously a
strong\ left ideal of $S$. By \cite[Lemma 1.1]{Ili2012}, there exists a
congruence $\sigma _{M}\in \mathrm{Cong}(M)$ on $M$ defined by%
\begin{equation*}
m\,\sigma _{M}\,m^{\prime }\Longleftrightarrow (0:_{S}m)=(0:_{S}m^{\prime })%
\text{.}
\end{equation*}%
Since $_{S}M$ is simple, $\sigma _{M}=$ $\vartriangle _{M}$ and, since $%
\{0\} $ and $S$ are the only strong left ideals of $S$, the annihilator $%
(0:_{S}m)=0$ for every $0\neq m\in M$ and, therefore, $M=\{0,m\}$. Let $N\in
|_{S}\mathcal{M}|$ be an essential extension of $M$ with the canonical
injection $i:M\rightarrowtail N$, and consider the congruence $\diamond
_{N}\in \mathrm{Cong}(N)$ on $N$ defined by%
\begin{equation*}
n_{1}\diamond _{N}n_{2}\Longleftrightarrow l_{1}n_{1}=n_{2}+n_{2}^{\prime }%
\text{ }\&\text{ }l_{2}n_{2}=n_{1}+n_{1}^{\prime }\text{ for some }%
l_{1},l_{2}\in \mathbb{N},\text{ }n_{1}^{\prime },n_{2}^{\prime }\in N\text{.%
}
\end{equation*}%
By \cite[Lemma 2.2]{Ili2012}, $N^{\diamond }:=N/\diamond _{N}$ is an
additively idempotent semimodule with the canonical surjection $%
p:N\twoheadrightarrow N^{\diamond }$. As $m\in I^{+}(M)$, it is clear that $%
(m,0)\notin \diamond _{N}$ and therefore the map
$pi:M\rightarrowtail N\twoheadrightarrow N^{\diamond }$ is
injective. Since $N$ is an essential extension of $M$, one has that
$p$ is an injective surjection, \textit{i.e.}, $p$ is an isomorphism
and $N$ is an additively idempotent semimodule. Then, considering
the congruence $\sigma_{N}\in \mathrm{Cong}(N)\ $on $N$
and using the same arguments as above for the semimodule $M$, one sees that $%
N/\sigma _{N}=\{\overline{0},\overline{m}\}$ with the natural surjection $%
\pi :N\longrightarrow N/\sigma _{N}$. Noting that the map $\pi
i:M\longrightarrow N\longrightarrow N/\sigma _{N}$ is injective and $N$ is
an essential extension of $M$, one concludes that $\pi $ is an isomorphism, $%
M=N$, and, applying Theorem 3.2 (2), ends the proof.\textit{\ \ \ \ \ \ }$%
_{\square }\medskip $

It is easy to see that $V(S)$ is a strong left and right ideal in a semiring
$S$. From this observation, Theorem 3.1 (3) and Proposition 3.3, we
have\medskip

\noindent \textbf{Corollary 3.4}\label{zeroic-cor} \textit{A semiring with
only two strong left (right) ideals is a left (right) V-semiring iff it is
either a left (right) V-ring, or a zeroic proper semiring.\medskip }

It is obvious that in any proper division semiring there are only two
trivial strong left (right) ideals, and therefore from Corollary 3.4 we
obtain\medskip

\noindent \textbf{Corollary 3.5}\label{ds-V} \textit{A division semiring}
\textit{is a left (right) V-semiring iff it} \textit{is either a zeroic
proper division semiring or a division ring.\medskip }

Now, we introduce a quite interesting class of semirings, naturally
extending the class of all rings. For any semiring $S$, let $%
P(S):=V(S)\,\cup \,\{1+s\,|\,s\in S\}$. It is easy to see that $P(S)$ is
always a subsemiring of $S$; and when $P(S)=S$, we say that the semiring $S$
is \emph{anti-bounded}.

From Proposition 3.3 we immediately have our first observation about
anti-bounded semirings.\medskip

\noindent \textbf{Corollary 3.6}\label{z-anti-b-V} \textit{A zerosumfree
anti-bounded semiring} $S$ \textit{is a left (right) V-semiring iff it}
\textit{is zeroic.\medskip }

\noindent \textit{Proof} Indeed, a zerosumfree anti-bounded semiring $S$ has
only two strong left (right) ideals: If an ideal $I\neq 0$ is a strong left
(right) ideal of $S$, $1+s\in I$ for some $s\in S$ and, hence, $1\in I$.%
\textit{\ \ \ \ \ \ }$_{\square }\medskip $

Clearly, proper division semirings are zerosumfree semirings containing only
two strong left (right) ideals. However, the class of zerosumfree semirings
possessing only two strong left (right) ideals, as the following example
shows, is quite wider than the class of proper division semirings.\medskip

\noindent \textbf{Example 3.7}\label{Bn} Let $n$ be a nonzero natural number
and $\mathbf{B}_{n+1}$ the join-semilattice defined on the chain $0<1<\cdots
<n$. Equip $\mathbf{B}_{n+1}$ with a structure of a semiring with addition $%
x+y:=x\vee y$ and multiplication%
\begin{equation*}
xy:=\left\{
\begin{array}{lcl}
0\text{,} &  & \text{if }x=0\,\,\text{or}\,\,y=0 \\
&  &  \\
x\vee y\text{,} &  & \text{otherwise}%
\end{array}%
\right. \text{.}
\end{equation*}%
It is easy to see that $\mathbf{B}_{n+1}$ is a zerosumfree anti-bounded
zeroic semiring with infinite element\textit{\ }$\infty =n$,\textit{\ }that
also is, by Corollary 3.6, a V-semiring and is not a division semiring. Of
course, $\mathbf{B}_{2}$ coincides with Boolean semiring $\mathbf{B}$%
.\medskip

Recall (see, for example, \cite{Kat2004} and \cite{KN2011}) that two
semirings $T$ and $S$ are said to be \textit{Morita equivalent} if the
semimodule categories $_{T}\mathcal{M}$ and $_{S}\mathcal{M}$ are equivalent
categories; \textit{i.e.}, there exist two additive functors $F:$ $_{T}%
\mathcal{M}\longrightarrow $ $_{S}\mathcal{M}$ and $G:$ $_{S}\mathcal{M}%
\longrightarrow $ $_{T}\mathcal{M}$, and natural isomorphisms $\eta
:GF\longrightarrow Id_{_{T}\mathcal{M}}$ and $\xi :FG\longrightarrow Id_{_{S}%
\mathcal{M}}$. By \cite[Theorem 4.12]{KN2011}, two semirings $T$ and $S$ are
Morita equivalent iff the semimodule categories $\mathcal{M}_{T}$ and $%
\mathcal{M}_{S}$ are equivalent categories. Following \cite{KN2011}, a left
semimodule $_{T}P\in |_{T}\mathcal{M}|$ is said to be a \textit{generator}
in the category $_{T}\mathcal{M}$ if the regular semimodule $_{T}T\in |_{T}%
\mathcal{M}|$ is a retract of a finite direct sum $\oplus _{i}P$ of the
semimodule $_{T}P$; and a left semimodule $_{T}P\in |_{T}\mathcal{M}|$ is
said to be a \textit{progenerator} in $_{T}\mathcal{M}$ if it is a finitely
generated projective generator. Finally, by \cite[Theorem 4.12]{KN2011} two
semirings $T$ and $S$ are Morita equivalent iff there exists a progenerator $%
_{T}P\in |_{T}\mathcal{M}|$ in $_{T}\mathcal{M}$ such that the
semirings $S$ and $\mathrm{End}(_{T}P)$ are isomorphic.

Our next observation is that the classes of simple and injective semimodules
are preserved by Morita equivalences of semirings, more precisely:\medskip

\noindent \textbf{Lemma 3.8} \textit{Let} $F:$ $_{T}\mathcal{M}%
\rightleftarrows $ $_{S}\mathcal{M}:G$ \textit{be an equivalence of the
semimodule categories} $_{T}\mathcal{M}$ \textit{and} $_{S}\mathcal{M}$%
\textit{. Then, a semimodule} $M\in |_{T}\mathcal{M}|$ \textit{is simple
(injective) iff} $F(M)\in |_{S}\mathcal{M}|$ \textit{is simple
(injective).\medskip }

\noindent \textit{Proof} By the dual of \cite[Lemmas 4.7 and 4.10]{KN2011},
one easily has that a semimodule $M\in |_{T}\mathcal{M}|$ is injective if
and only if $F(M)\in |_{S}\mathcal{M}|$ is injective.

Now, let $M\in |_{T}\mathcal{M}|$ be a simple semimodule, and $%
f:F(M)\longrightarrow N$ be a nonzero homomorphism in $_{S}\mathcal{M}$.
Then, $G(f):M\cong G(F(M))\longrightarrow G(N)$ is a nonzero homomorphism in
$_{T}\mathcal{M}$, and, hence, it is injective. Applying the functor $F$ and
the dual of \cite[Lemmas 4.7 ]{KN2011}, we have
\begin{equation*}
F(G(f)):F(M)\cong F(G(F(M)))\longrightarrow F(G(N))\cong N
\end{equation*}
is injective, and, hence, $f$ is injective, too; and by Proposition 2.5, $%
F(M)$ is a simple semimodule, too.\textit{\ \ \ \ \ \ }$_{\square }\medskip $

Applying Lemma 3.8, we immediately establish that for semirings to be a left
(right) V-semiring is a Morita invariant property.\medskip

\noindent \textbf{Theorem 3.9}\label{Mor} \textit{Let }$T$\textit{\ and }$S$%
\textit{\ be Morita equivalent semirings. Then }$T$\textit{\ and }$S$\textit{%
\ are left (right) V-semiring simultaneously.\medskip }

Following \cite{EHJK2001}, we call a semiring $S$ \emph{congruence-simple}%
\textit{\ }(\emph{ideal-simple}) if the diagonal and universal congruences
are the only congruences on $S$ (if $0$ and $S$ are the only ideals of $S$);
finally, a semiring $S$ is said to be \emph{simple} if it is simultaneously
congruence- and ideal-simple. In contrast to the varieties of groups and
rings, the research on simple semirings has just recently been started, and
therefore not too much is so far known on them (for some recent activity and
results on this subject, one may consult \cite{EHJK2001}, \cite{Mon2004},
\cite{EK2007}, \cite{Zum2008}, \cite{JKM2009}, \cite{KNT2011}, \cite{KZ2013}%
). Thus, regarding relations between V-semirings and different variations of
`simplicity' of semirings, the following observations deserve to be
mentioned.\medskip

\noindent \textbf{Corollary 3.10} \textit{A simple semiring} $S$ \textit{%
with either an infinite element or a projective minimal one-sided ideal is a
left (right) V-semiring.}\medskip

\noindent \textit{Proof} \ For a semiring $S$ with an infinite element, the
result immediately follows from Example 3.7, Theorem 3.9 and \cite[Theorem
5.7]{KNZ}.

If a semiring $S$ possesses a projective minimal one-sided ideal, then the
result follows from Example 3.7, Theorem 3.9 and \cite[Theorems 5.7 and 5.11]%
{KNZ}. \textit{\ \ \ \ \ \ }$_{\square }\medskip $

Using the complete description of ideal-simple and congruence-simple
Artinian subtractive semirings given in \cite[Theorems 3.7 and 4.5]{KNT2011}%
, one has\medskip

\noindent \textbf{Corollary 3.11} \textit{An ideal-simple left (right)
Artinian, left (right) subtractive semiring }$S$\textit{\ is a left (right)
V-semiring iff }$S\cong M_{n}(D)$\textit{\ for some division ring }$D$%
\textit{, or }$S$\textit{\ is a zeroic division semiring. Any
congruence-simple left (right) Artinian,  left (right) subtractive semiring
is a  left (right) V-semiring.\medskip }

\noindent \textit{Proof} \ The first statement follows right away from \cite[%
Theorem 3.7]{KNT2011}, Corollary 3.5, Theorem 3.9, \cite[Theorem 5.14]%
{Kat2004}, and the classical fact --- over division rings all modules are
injective.

The second statement also follows right away from \cite[Theorem 4.5]{KNT2011}%
, Example 3.7, \cite[Theorem 5.14]{Kat2004}, and the same classical fact.
\textit{\ \ \ \ \ \ }$_{\square }\medskip $

As usual, a semiring $S$ is said to be \emph{left }(\emph{right}) \emph{%
semisimple} if the regular semimodule is a direct sum of left (right) atom
ideals. Recall (see, for example, \cite[Theorem 7.8]{HW1996}, or \cite[%
Theorem 4.5]{KNT2009}) that a semiring $S$\ is (left, right) \textit{%
semisimple} iff
\[
S\cong M_{n_{1}}(D_{1})\times \cdots \times M_{n_{r}}(D_{r})\text{,}
\eqno(*)
\]%
where $M_{n_{1}}(D_{1}),\ldots ,M_{n_{r}}(D_{r})$ are the semirings of $%
n_{1}\times n_{1},\ldots ,n_{r}\times n_{r}$\ matrices for suitable division
semirings $D_{1},\ldots ,D_{r}$\ and positive integers $n_{1},\ldots ,n_{r}$%
. In the sequel, we will refer to such an isomorphism $(\ast )$ as a \textit{%
direct product representation} of a semisimple semiring $S$.

Our next result provides us with a description of semisimple left
(right)\linebreak V-semirings:\medskip

\noindent \textbf{Theorem 3.12}\label{semisimple-V} \textit{The following
conditions for a semisimple semiring }$S$\textit{\ are equivalent:}

\textit{(1) }$S$\textit{\ is a left (right) V-semiring;}

\textit{(2) }$S\cong M_{n_{1}}(D_{1})\times \cdots \times M_{n_{r}}(D_{r})$%
\textit{,\ where }$D_{1},\cdots ,D_{r}$\textit{\ are either division rings
or zeroic division semirings.\medskip }

\noindent \textit{Proof} \ (1) $\Longrightarrow $ (2). Let $S$ be a left
V-semiring and $S\simeq M_{n_{1}}(D_{1})\times \cdots \times
M_{n_{r}}(D_{r}) $ a direct product representation of $S$. By Theorem 3.1, $%
M_{n_{i}}(D_{i})$, $i=1,\cdots ,r$, are left V-semirings, too. Whence, by
\cite[Theorem 5.14]{Kat2004} and Theorem 3.9, $D_{i}$, $i=1,\cdots ,r$, are
also left (right) V-semiring, and applying Corollary 3.5 one gets the
statement.

(2) $\Longrightarrow $ (1). This implication follows straight away from
Corollaries 3.2, 3.5, \cite[Theorem 5.14]{Kat2004} and Theorem 3.9. \textit{%
\ \ \ \ \ \ }$_{\square }\medskip $

\noindent \textbf{Corollary 3.13}\label{ad-id-semi} \textit{Every additively
regular, in particular every finite, semisimple semiring is a left \emph{(}%
right\emph{)} V-semiring.\medskip }

\noindent \textit{Proof} \ First, let $S\simeq M_{n_{1}}(D_{1})\times \cdots
\times M_{n_{r}}(D_{r})$ be a direct product representation of an additively
regular semisimple semiring $S$. By Corollary 3.2, it is easy to see that,
without loss of generality, we can consider only the case when all $D_{i},$ $%
i=1,\cdots ,r$, are proper additively regular and, therefore, even
additively idempotent (see also \cite[p. 4349]{KNT2011}), division
semirings. From the latter, Corollary 3.5, Theorem 3.9, and \cite[Theorem
5.14]{Kat2004} one gets that all semirings $M_{n_{i}}(D_{i}),$ $i=1,\cdots ,r
$, are left (right) V-semirings, too, and our assertion follows from
Corollary 3.2.

If $S$ is a finite semisimple semiring, then all division semirings $D_{i},$
$i=1,\cdots ,r$, are also finite. Then, applying, for instance, \cite[%
Proposition 1.2.3]{How1995}, one gets that the monogenic additive
subsemigroup $<1_{i}>$ of each monoid $(D_{i},+,0)$ generated by $1_{i}\in
D_{i},$ $i=1,\cdots ,r$, contains a nonzero idempotent, and therefore, each
division semiring $D_{i},$ $i=1,\cdots ,r$, is additively idempotent. So, $S$
is also an additively idempotent semiring, and the statement follows from
the previous one above. \textit{\ \ \ \ \ \ }$_{\square }\medskip $

We conclude this section by considering some relations between V-semirings
and the \emph{Jacobson-Bourne radical for semirings} --- a semiring analog
and/or generalization of the classical Jacobson radical for rings ---
introduced by S. Bourne in \cite{Bou1951}.

Recall \cite{Bou1951} that a right ideal $I$ of a semiring $S$ is said to be
\emph{right semiregular} if, for every pair of elements $i_{1},i_{2}$ in $I$%
, there exist elements $j_{1}$ and $j_{2}$ in $I$ such that $%
i_{1}+j_{1}+i_{1}j_{1}+i_{2}j_{2}=i_{2}+j_{2}+i_{1}j_{2}+i_{2}j_{1}$. As was
shown in \cite[Theorems 3 and 4]{Bou1951}, the sum of all right semiregular
ideals of a semiring $S$, denoted by $\mathrm{J}(S)$, is also a right
semiregular ideal of $S$. This ideal is called the \emph{Jacobson radical}
of the semiring $S$, and $S$ is said to be \emph{Jacobson-semisimple} if $%
\mathrm{J}(S)=0$. As was shown in \cite[Theorem 3.75]{Lam1999}, all left
(right) V-ring $S$ are Jacobson-semisimple rings, \textit{i.e.}, $\mathrm{J}%
(S)=0$. However, this is not true for V-semirings in general: For example,
by Corollary 3.5, the Boolean semifield $\mathbf{B}$ is a left (right)
V-semiring, but it is very easy to see that \textrm{J}$(\mathbf{B})=\mathbf{B%
}$ (see also \cite[Example 3.7]{KN}). In light of this fact, our next result
shows that the class of the Jacobson-semisimple V-semirings contains only
rings, namely:\medskip

\noindent \textbf{Theorem 3.14} \textit{A left (right) V-semiring} $S$
\textit{is Jacobson-semisimple iff }$S$\textit{\ is a left (right)
V-ring.\medskip }

\noindent \textit{Proof} $\Longrightarrow $. Let $S$ be a
Jacobson-semisimple left (right) V-semiring. Then, by
\cite[Corollary 4.6]{KN}, $S$ is semiisomorphic to a subdirect
product of some additively cancellative semirings
$\{S_\lambda\}_{\lambda \in \Lambda}$ whose rings of differences
$S_\lambda^\Delta$ (see, \cite[p. 101]{Gol1999}) are
isomorphic to dense subrings of the rings of linear transformations $\mathrm{%
End}(_{D_{\lambda }}V_{\lambda })$ of vector spaces $_{D_{\lambda
}}V_{\lambda }$ over division rings $D_{\lambda }$ for each $\lambda
\in \Lambda $. Since all semirings $S_{\lambda }$, $\lambda \in
\Lambda $, are additively cancellative, by Theorem~3.1 (4), they are
clearly left (right) V-rings, and, hence, $S$ is a left (right)
V-ring, too.

$\Longleftarrow $. This follows from \cite[Theorem 3.75]{Lam1999}. \textit{\
\ \ \ \ \ }$_{\square }\medskip $

There is another very natural semiring analog of the Jacobson radical for
semirings, $\mathrm{J}_{s}$-radical, based on the class of simple
semimodules, considered in \cite{KN}, and is, in general, different from $%
\mathrm{J}$-radical (\cite[Example 3.7]{KN}). It is easy to see that a
proper division semiring $D$ is a left (right) V-semirings with $\mathrm{J}%
_{s}(D)=0$, and hence, an analog of Theorem 3.14 for the $\mathrm{J}_{s}$%
-semisimple semirings is not true. In light of this, we end this section by
posting, in our view, an interesting and perspective\medskip\ question.

\noindent \textbf{Problem 1.} Describe all $\mathrm{J}_{s}$-semisimple
V-semirings.

\section{Semirings all of whose cyclic semimodules are injective}

\qquad Inspired by the well-known characterization of semisimple rings as
the rings all of whose cyclic modules are injective given by Barbara Osofsky
(\cite{Oso1964}, or \cite[Corollary 6.47)]{Lam1999}), we initiate in this
section a study of the so called \emph{left }(\emph{right}) \emph{%
CI-semirings }--- semirings all of whose cyclic semimodules are injective.
As our next observation shows, the class of CI-semirings is significantly
wider than that of semisimple rings; and the CI-semirings, in our view, will
constitute a very interesting direction of future investigations in the
non-additive semiring setting. Here, by characterizing the CI-semirings
within some special and important classes of semirings, we are just
starting, and hopefully motivating too, the research in this promising
direction.\medskip

\noindent \textbf{Proposition 4.1} \textit{A semiring} $S$ \textit{is a left
(right) CI-semiring iff} $S=R\oplus T$\textit{\ with} $R$\textit{\ and} $T$
\textit{a semisimple ring and a zerosumfree left (right) CI-semiring with an
infinite element, respectively.\medskip }

\noindent \textit{Proof} $\Longrightarrow $. Let $S$ be a left CI-semiring, $%
\equiv _{V(S)}$ the Bourne congruence on $S$. It is clear that the factor
semiring $\overline{S}:=S/\equiv _{V(S)}$ is zerosumfree, and the natural
surjection $\pi :S\longrightarrow \overline{S}$ induces the \emph{%
restriction functor }$\pi ^{\#}:$ $_{\overline{S}}\mathcal{M}\longrightarrow
$ $_{S}\mathcal{M}$ (see \cite[p. 202, Proposition 4.1]{Kat2004}). The
semiring $\overline{S}$ is a left CI-semiring too: Indeed, if $M$ $\in |_{%
\overline{S}}\mathcal{M}|$ is a cyclic $\overline{S}$-semimodule, then $\pi
^{\#}(M)\in |_{S}\mathcal{M}|$ is a cyclic $S$-semimodule as well and,
hence, injective; for \cite[Lemma 5.2]{KN2011}, the latter implies that $M$
is an injective $\overline{S}$-semimodule as well. In particular, for the
zerosumfree semiring $\overline{S}$, the regular semimodule $_{\overline{S}}%
\overline{S}\in |_{\overline{S}}\mathcal{M}|$ is injective, and therefore,
by \cite[Proposition 1.7]{Ili2010}, the semiring $\overline{S}$ contains an
infinite element and, hence, it is a zeroic zerosumfree semiring. So, by
\cite[Proposition 2.9]{Ili2012}, $S=R\oplus T$, where $R$ is a ring and $T$
is a semiring isomorphic to $\overline{S}$. One readily sees that $T$ is a
left CI-semiring with an infinite element and the ring $R$ is a left
CI-ring. Therefore, by \cite{Oso1964} (see also \cite[Theorem 1.2.9]{Lam2001}
or \cite[Corollary 6.47)]{Lam1999}), $R$ is a semisimple ring.

$\Longleftarrow $. It is easy to see that a finite direct product of left
CI-semirings is also a left CI-semiring. Using this observation and
Osofsky's result \cite{Oso1964}, one ends the proof. \textit{\ \ \ \ \ \ }$%
_{\square }\medskip $

From this proposition one sees right away that the problem of describing
CI-semirings is\ actually reduced to the quite nontrivial problem of
describing all zerosumfree left (right) CI-semirings with infinite elements.
One important subclass of the class of zerosumfree left (right) CI-semirings
with infinite elements --- the class of distributive lattices with $0$ and $1
$ (\textit{i.e.}, bounded distributive lattices) that are CI-semirings ---
have been already considered in Dr. Kornienko's PhD thesis, supervised by
late Prof. L. A. Skornjakov and successfully defended at Moscow State
University (MGU), USSR, in 1979 (see also the announcement of the results in
\cite[p. 118]{Kor1978}). However, since a proof of Dr. Kornienko's result,
to the best of our knowledge, has never been published in any publicly
available publications, we find it reasonable to present our independent
proof\ of that result here, too. First, we need the following easy
fact:\medskip

\noindent \textbf{Lemma 4.2 }(\textit{cf.} \cite[Theorem 32]{Sal1988})
\textit{Any infinite Boolean algebra }$B$\textit{\ has a countable set of
orthogonal idempotents.\medskip }

\noindent \textit{Proof} \ Let $B=(B;+,\cdot ,^{\prime },0,1)$ be an
infinite Boolean algebra. Take an element $a\in B$ such that $0\neq a\neq 1$%
; obviously, $0\neq a^{\prime }\neq 1$ and $B=aB\oplus a^{\prime }B$. It is
clear that $aB\ $and $a^{\prime }B$ are also Boolean algebra and at least
one of them, say $B_{1}=a^{\prime }B$, is infinite. Let $e_{1}:=a$.
Repeating the same procedure for the algebra $B_{1}=e_{1}^{\prime }B$, one
gets a nonzero element $e_{2}\in $ $B_{1}$ and an infinite subalgebra $%
B_{2}=e_{2}^{\prime }B_{1}$. Continuing this process, one readily gets a
countable sequence $\{e_{1},e_{2},\ldots \}$ of mutually orthogonal
idempotents of $B$. \textit{\ \ \ \ \ \ }$_{\square }\medskip $

The next result, announced in \cite[p. 118]{Kor1978}, characterizes
the CI-semi\-rings within the class of bounded distributive
lattices.\medskip

\noindent \textbf{Theorem 4.3} \textit{A bounded distributive lattice is a
CI-semiring iff it is a finite Boolean algebra.\medskip }

\noindent \textit{Proof }$\Longrightarrow $. Let a CI-semiring $B$ be a
bounded distributive lattice. By \cite[Theorem]{Kor1979}, the lattice $B$ is
a complete Boolean algebra. It is clear that a factor algebra of a Boolean
algebra as well as a factor semiring of \ a CI-semiring are a Boolean
algebra and a CI-semiring, respectively, too. Therefore, it is enough to
show that any infinite complete Boolean algebra has a non-complete factor
algebra. Thus, assume that the lattice $B$ is an infinite complete Boolean
algebra with, by Lemma 4.2, a countable sequence $\{e_{1},e_{2},\ldots \}$
of mutually orthogonal nonzero idempotents. We shall show that the factor
algebra $\overline{B}:=B/A$ of the algebra $B$ with respect to the ideal $%
A:=\Sigma _{i=1}^{\infty }e_{i}B$ is not a complete Boolean algebra.

Consider $I_{p}:=\{$ $p^{k}|$ $k=1,2,\ldots \}$ and $f_{p}:=\vee _{i\in
I_{p}}e_{i}\in B$ for every prime number $p$. Certainly, the set $\overline{F%
}=\{$ $\overline{f_{p}}$ $|$ $p$ is prime$\}$ of the elements of the Boolean
algebra $\overline{B}$ has an upper bound (for instance, the element $%
\overline{1}$ is its upper bound); however, as we will show, there is no
supremum of $\overline{F}$ --- for every upper bound $\overline{f}$ of $%
\overline{F}$, there exists another upper bound $\overline{g}$ of the set $%
\overline{F}$ such that $\overline{g}<$ $\overline{f}$. Indeed, from the
inequality $\overline{f_{p}}\leq \overline{f}$ for every prime $p$, it
follows that there are some elements $a_{p1},a_{p2}\in A$ such that $\vee
_{i\in I_{p}}e_{i}+f+a_{p1}=f_{p}+f+a_{p1}=f+a_{p2}$. Whence, multiplying
the latter equation by $e_{i}$, $i\in I_{p}$, and using the join infinite
distributive identity for complete Boolean algebra \cite[Lemma 2.4.10]%
{Gra1998} and the mutual orthogonality of the idempotents $e_{i}$, $%
i=1,2,\ldots $, we have $e_{i}=fe_{i}+a_{p2}e_{i}$ and $a_{p2}e_{i}=0$ for
all but finite $i\in I_{p}$. From the latter, $e_{i}=fe_{i}$ for all but
finite $i\in I_{p}$; therefore, for the index set $J_{p}:=\{$ $i\in I_{p}|$ $%
e_{i}=fe_{i}$ $\}$ we have $|I_{p}\backslash $ $J_{p}|<\infty $ and $%
\overline{f_{p}}=\overline{g_{p}}$ for $g_{p}:=\vee _{i\in J_{p}}e_{i}$.

Now, let $g:=\vee _{_{p}}g_{p}$. Since $\overline{f_{p}}=\overline{g_{p}}%
\leq \overline{g}$ for each $p$, $\overline{g}$ is an upper bound of the set
$\overline{F}$ and, by construction, $\overline{g}\leq $ $\overline{f}$. If
the latter inequality is a strong one, $\overline{g}$ is exactly the element
we are looking for. So, suppose that $\overline{g}=$ $\overline{f}$, and for
each prime $p$ consider the index sets $K_{p}:=J_{p}\backslash $ $\{j_{p}\}$%
, where $j_{p}$ is the smallest element of the index set $J_{p}$. For $%
|I_{p}\backslash K_{p}|<\infty$, we have
$\overline{f_p}=\overline{h_p}$ for $h_p:=\vee _{i\in K_p}e_i$ and
$\overline{f_p}=\overline{h_p}\leq\overline{h}$ for $h:=\vee
_{_{p}}h_p$; whence, $\overline{h}$\ is an upper bound of
$\overline{F}$. Moreover, $g=h+u$, where $u:=\vee _{_{p}}e_{j_{p}}$
and, applying again \cite[Lemma 2.4.10]{Gra1998}, $hu=0$; therefore,
$\overline{g}=\overline{h}+\overline{u}$ and
$\overline{h}\overline{u}=\overline{0}$. Whence,
$\overline{h}<\overline{g}$ provided $\overline{u}\neq
\overline{0}$. However, the latter is almost obvious: Indeed, there
are infinite number of prime numbers $p$, all $j_{p}$ are different
for different prime numbers; whence, $u\notin A$ and
$\overline{u}\neq \overline{0}$.

$\Longleftarrow $. This immediately follows from \cite[Theorem 4]{Fof1972}.
\textit{\ \ \ \ \ \ }$_{\square }\medskip $

Combining \cite[Corollary 8]{Kor1979}, \cite[Theorem]{Kor1980}, and
Theorem 4.3, one obtains the following characterization of finite
Boolean algebras among bounded distributive lattices:\medskip

\noindent \textbf{Corollary 4.4} \textit{The following conditions for a
bounded distributive lattice }$B$\textit{\ are equivalent:}

\textit{(1) All cyclic }$B$\textit{-semimodules are projective;}

\textit{(2) All atom }$B$\textit{-semimodules are projective;}

\textit{(3) All subsemimodules of the regular semimodule }$B_{B}$\textit{\
are injective;}

\textit{(4) }$B$\textit{\ is a CI-semiring;}

\textit{(5) }$B$\textit{\ is a finite Boolean algebra.\medskip }

As usual, we denote by $U(S)$ the set of all units of a semiring $S$ and,
following \cite[p. 56]{Gol1999}, we say that a semiring $S$ is a \emph{%
Gelfand semiring} if the element $1+s\in $ $U(S)$ for every element $s\in S$%
. Of course, bounded distributive lattices are among Gelfand semirings. But
the class of the Gelfand semirings is quite wider as the following example,
for instance, shows:\medskip

\noindent \textbf{Example 4.5} \cite[Example 4.49]{Gol1999} Let $A$ be a
nonempty set having more than one element and $S=(\mathbb{R^{+}})^{A}$ the
set of all functions from $A$ to $\mathbb{R}^{+}$, where $\mathbb{R}^{+}$ is
the set of all nonnegative real numbers, with the canonical semiring
structure defined on it. Obviously, $U(S)=\{f\in S\mid $ $f(a)>0$ for all $%
a\in A\}$ and $S$ is a Gelfand semiring that is not a bounded distributive
lattice.\medskip

In this connection, the following characterization of CI-, Gelfand semirings
is certainly of interest. \medskip

\noindent \textbf{Theorem 4.6} \textit{A Gelfand semiring} $S$ \textit{is a
left (right) CI-semiring iff it is a finite Boolean algebra.\medskip }

\noindent \textit{Proof }$\Longrightarrow $. Let $S$ be a left CI-semiring.
By Proposition 4.1, $S=R\oplus T$ with $R$ a semisimple ring and $T$ a left
CI-semiring with an infinite element $\infty $. Let $1=r+t$ for some $r\in R$
and $t\in T$. Since $S$ is a Gelfand semiring, $t=1-r\in U(S)$ and, by \cite[%
Proposition 4.50]{Gol1999}, $\infty =t+\infty \in U(S)$. As $\infty +\infty
^{2}=\infty $, we have $\infty =1+\infty =1$, $\ $and therefore $S=T$ is an
additively idempotent semiring.

Consequently, there is a natural partial order on the semiring $S$ given by $%
a\leq b$ iff $a+b=b$ with respect to which $a\vee b=a+b$ for all $a,b\in S$.
Since $S$ is a CI-semiring, by \cite[Theorem 2.2]{ASW1998}, $S$ is also a
von Neumann regular semiring, and hence, for any $s\in S$ there exists $x\in
S$ such that $s=sxs$. The latter implies that $%
s+s^{2}=sxs+s^{2}=s(x+1)s=s^{2}$ and $s\leq s^{2}$. On the other hand, as $%
s\leq 1$ one has $s^{2}\leq s\ $and therefore $s=s^{2}$ for all $s\in S$.
Moreover, $S$ is a commutative semiring: Indeed, since $a\leq 1$ and $b\leq
1\ $for all $a,b\in S$, it follows that $ab\leq b,$ $ab\leq a$ and $%
ab=(ab)^{2}=(ab)(ab)\leq ba$, and by symmetry, $ba\leq ab$ and, hence, $%
ab=ba $. We show that $a\wedge b=ab$ for all $a,b\in S$: First, $%
a+ab=a(1+b)=a\cdot 1=a$ and $b+ab=(1+a)b=1\cdot b=b$, and hence, $ab\leq a$,
$ab\leq b$, and $ab\leq a\wedge b$. Next, if for some $x\in S$ we have $%
x\leq a$ and $x\leq b$, then%
\begin{equation*}
\begin{tabular}{lllll}
$ab$ & $=$ & $(x+a)(x+b)$ & $=$ & $x^{2}+xb+ax+ab$ \\
& $=$ & $x+xb+xa+ab$ & $=$ & $x(1+b+a)+ab$ \\
& $=$ & $x\cdot 1+ab$ & $=$ & $x+ab$,%
\end{tabular}%
\end{equation*}%
\textit{i.e.}, $x\leq ab$, and therefore, $a\wedge b=ab$. Thus, our semiring
$S$ is, in fact, a bounded distributive lattice $(S,\vee ,\wedge )$, and
applying Theorem 4.3 one gets the statement.

$\Longleftarrow $. This is true by Theorem 4.3. \textit{\ \ \ \ \ \ }$%
_{\square }\medskip $

In light of these results, the following problem, we believe, seems to be
quite natural and interesting.\medskip

\noindent \textbf{Problem 2.} To which extent can the equivalent conditions
of Corollary 4.4 be extended to Gelfand semirings?\medskip

Also, in the connection with a description of the bounded distributive
lattices that are V-semirings obtained in \cite[Corollary 3.11]{Ili2012}, it
is natural to post\medskip

\noindent \textbf{Problem 3.} Describe Gelfand V-semirings.\medskip\

Next, using Theorem 4.6, we give a complete description of left subtractive
left CI-semirings. \medskip

\textbf{Theorem 4.7} \textit{A left subtractive semiring} $S$ \textit{is a
left CI-semiring iff} $S = R\oplus T$ \textit{with} $R$ \textit{and} $T$
\textit{a semisimple ring and a finite Boolean algebra, respectively.\medskip%
}

\noindent \textit{Proof } $\Longrightarrow$. By Proposition 4.1, $%
R=S\oplus T,$ where $S$ is a classical semisimple ring, and $T$ is a
zerosumfree left CI-semiring with the infinite element $\infty $. It is easy
to see (\textit{e.g.,} by \cite[Lemma 4.7]{KNT2009}) that $T$ is a left
subtractive semiring, too. Then, the left ideal $T\infty $ is subtractive in
$T$ and it follows from $1_{T}+\infty =\infty $ that $1_{T}\in T\infty $,
that means, $t\infty =1_{T}$ for some $t\in T$.

On the other hand, we have $\infty^2 +\infty =\infty$. This implies that
\begin{equation*}
\infty =\infty + 1_T = 1_T,
\end{equation*}
therefore, $T$ is a Gelfand semiring. By Theorem 4.6, $T$ is a finite
Boolean algebra.

$\Longleftarrow$. The semirings $R$ and $T$ are, obviously, left
subtractive semirings. By \cite[Lemma 4.7]{KNT2009}, $S=R\oplus T$
is a left subtractive semiring, too.

Applying Proposition 4.1 and Theorem 4.6, we conclude that $S$ is a left
CI-semiring.\textit{\ \ \ \ \ \ }$_{\square }\medskip $

By \cite[Theorem 14.1]{EHJK2001} and \cite[Corollary 4.6]{KNT2011},  the
concepts of congruence-simpleness and ideal-simpleness coincide for finite
commutative and finite left (right) subtractive semirings. Also, using
Theorem 4.6 and \cite[Theorem 3.7]{KNT2011}, we have the same situation for
subtractive CI-semirings  and, in fact, solve \cite[Problem 3]{KNT2011} and
\cite[Problem 5]{KNZ} in the class of subtractive CI-semirings.\medskip

\textbf{Corollary 4.8.} \textit{For a left subtractive semiring} $S$\textit{,%
} \textit{the following conditions are equivalent:}

\textit{(i)} $S$ \textit{is a congruence-simple left CI-semiring;}

\textit{(ii)} $S$ \textit{is an ideal-simple left CI-semiring;}

\textit{(iii)} $S\cong M_n(D)$ \textit{for some division ring} $D,$ \textit{%
or} $S\cong \mathbf{B}.\medskip$

\noindent \textit{Proof } (i) $\Longrightarrow$ (ii). It follows immediately
from \cite[Proposition 4.4]{KNT2011}.

(ii) $\Longrightarrow $ (iii). Assume that $S$ is an ideal-simple left
CI-semiring. Applying Theorem 4.7, $S$ is either a semisimple ring, or a
finite Boolean algebra. If $S$ is a semisimple ring, then since $S$ is
ideal-simple, $S\cong M_{n}(D)$ for some division ring $D$ and $n\geq 1$.
Otherwise, $S$ is a finite Boolean algebra without proper nonzero ideals;
therefore, $S$ is just the Boolean semifield $\mathbf{B}$.

(iii) $\Longrightarrow$ (i). It follows immediately from Proposition 4.1,
Theorem 4.7, and \cite[Theorem 4.5]{KNT2011}. \textit{\ \ \ \ \ \ }$%
_{\square }\medskip $

Now, let us consider semisimple CI-semirings. Using the direct
product representation $(*)$ of semisimple semirings, our
considerations are naturally reduced to the ones of the matrix
CI-semirings over division semirings about which the following
observation is crucial.\medskip\

\noindent \textbf{Proposition 4.9} \textit{The matrix semiring }$S=M_{n}(D)$%
\textit{\ over a division semiring }$D$\textit{\ is a left (right)
CI-semiring iff }$D$\textit{\ is a division ring, or }$D\cong \mathbf{B}$%
\textit{\ and }$n=1,2$\textit{.\medskip }

\noindent \textit{Proof \ } We will use the equivalence of the semimodule
categories $_{S}\mathcal{M}$ and $_{D}\mathcal{M}$ established in \cite[%
Theorem 5.14]{Kat2004}: $F:$ $_{S}\mathcal{M}\rightleftarrows $ $_{D}%
\mathcal{M}:G$, $F(A)=E_{11}A$ and $G(B)=B^{n}$, where $E_{11}\ $is the
matrix unit in $M_{n}(D)$.

$\Longrightarrow $. Let $S=M_{n}(D)$ be a left CI-semiring. By \cite[Theorem
2.2]{ASW1998}, $S$ is a von Neumann regular semiring, and, if $n\geq 3$, by
\cite[Proposition 2]{Ili2001}, $D$ is a regular ring. Hence, we need to
consider only the case with $n=1,2$.

Let $D$ be a zerosumfree division semiring. Since $_{D}D\cong F(G(_{D}D))$, $%
G(_{D}D)=$ $_{S}D^{2}$, and noting that $_{S}D^{2}$ is obviously a cyclic $S$%
-semimodule and, therefore, an injective one, by Lemma 3.8, one has that $%
_{D}D$ is an injective $D$-semimodule. Whence, by \cite[Proposition 1.7]%
{Ili2010}, $D$ has an infinite element $\infty $ such that $\infty
^{2}+\infty =\infty $, hence, $1=\infty +1=\infty $. So, $d^{-1}+1=1$ and,
hence, $d=1+d=1$ for every $0\neq d\in D$; therefore, $D=\{0,1\}\cong
\mathbf{B}.$

$\Longleftarrow $ . If $D$ is a division ring, then $S$ is a semisimple ring
and the statement is obvious. Let $n=1,2$, $D\cong \mathbf{B}$, and $M $ be
a cyclic left $M_{2}(\mathbf{B})$-semimodule with an $M_{2}(\mathbf{B})$%
-surjection $f:M_{2}(\mathbf{B})\twoheadrightarrow M.$ By \cite[Lemma 4.7]%
{KN2011}, $F(f):\mathbf{B}^{2}\cong F(M_{2}(\mathbf{B}))\twoheadrightarrow
F(M)$ is a $\mathbf{B}$-surjection and, hence, $\mathbf{B}^{2}/\equiv
_{F(f)}\cong F(M)$. It is clear that $\{(0,0)\},$ $\mathbf{B},$ $%
\{(0,0),(1,0),(1,1)\}$ and $\mathbf{B}^{2}$ are, up to isomorphism, the only
quotient semimodules of $\mathbf{B}^{2}$ which are injective by \cite[%
Theorem 4]{Fof1972}. Whence, $F(M)$ is an injective left $\mathbf{B}$%
-semimodule too, and therefore, by Lemma 3.8, $M\cong G(F(M))$ is an
injective left $M_{2}(\mathbf{B})$-semimodule as well. \textit{\ \ \ \ \ \ }$%
_{\square }\medskip $

Using this proposition, we obtain a complete description of
semisimple \linebreak CI-semirings:\medskip

\noindent \textbf{Theorem 4.10}\label{CI-semisimple} \textit{A semisimple
semiring} $S$ \textit{is a left (right) CI-semiring iff }$S\cong S_{1}\times
\cdots \times S_{r}$\textit{,} \textit{where} \textit{each} $S_{i}$, $%
i=1,\cdots ,r$, \textit{is either an Artinian simple ring, or isomorphic to}
$M_{n}(\mathbf{B})$ \textit{with} $n=1,2$\textit{.\medskip }

\noindent \textit{Proof \ } $\Longrightarrow $. Let
$M_{n_1}(D_1)\times \cdots \times M_{n_r}(D_r)$ be a direct product
representation $(*)$ of $S$. By \cite[Lemma 5.2]{KN2011},
$M_{n_i}(D_i)$, $i=1,\cdots ,r$, are left CI-semirings, and, by
Proposition 4.9, each $D_i$ is either a division ring or $D_i\cong
\mathbf{B}$ with $n=1,2$.

$\Longleftarrow $. This implication follows right away from Proposition 4.9
and the obvious fact that a finite direct product of left (right)
CI-semirings is also a left (right) CI-semiring. \textit{\ \ \ \ \ \ }$%
_{\square }\medskip $

A complete description of anti-bounded CI-semirings constitutes our next
main goal in this paper. We need first to justify the following important
and useful facts about some semirings, in particular about the ones
introduced in Example 3.7 above.\medskip

\noindent \textbf{Fact 4.11}\label{R3} \textit{The semiring} $\mathbf{B}_{3}$%
\textit{, defined on the chain }$0<1<$\textit{\ }$2$\textit{\ in Example 3.7,%
} \textit{is a CI-semiring.\medskip }

\noindent \textit{Proof \ } First note that, up to isomorphism, there are
only two nonzero cyclic $\mathbf{B}_{3}$-semimodules, namely $\{0,2\}$ and $%
\mathbf{B}_{3}$. Indeed, let $M=\mathbf{B}_{3}m$ be a nonzero cyclic $%
\mathbf{B}_{3}$-semimodule for some $0_{M}\neq m\in M$. Then, $2m=0_{M}\ $%
implies $0_{M}=2m=(1+2)m=m+2m=m+0_{M}=m$ which contradicts $0_{M}\neq m$; if
$m=2m\ $or $m\neq 2m\neq 0$, it is easy to see that then$\ M\simeq \{0,2\}$
or $M\simeq \mathbf{B}_{3}$, respectively.

Now, considering the category of semimodules $_{\mathbf{B}_{3}}\mathcal{M}$
over the additively idempotent commutative semiring $\mathbf{B}_{3}$, noting
that the regular semimodule $_{\mathbf{B}_{3}}\mathbf{B}_{3}\in |_{\mathbf{B}%
_{3}}\mathcal{M}|$ is a free, and therefore, a flat semimodule, and applying
\cite[Theorem 3.9 and Proposition 4.1]{Kat1997}, one obtains that the
`character' semimodule $\mathbf{2}^{\mathbf{B}_{3}}:=\mathcal{M}(\mathbf{B}%
_{3},\mathbf{2)}\in |_{\mathbf{B}_{3}}\mathcal{M}|$ (where $\mathbf{2}$ is
the two-element semilattice) of \ the semimodule $_{\mathbf{B}_{3}}\mathbf{B}%
_{3}$ is an injective semimodule. We conclude the proof by noting that the
semimodules $\mathbf{2}^{\mathbf{B}_{3}}$ and $_{\mathbf{B}_{3}}\mathbf{B}%
_{3}$ are obviously isomorphic, and that the $\mathbf{B}_{3}$-semimodules $%
\{0,2\} $ is a retract of $_{\mathbf{B}_{3}}\mathbf{B}_{3}$. \textit{\ \ \ \
\ \ }$_{\square }$\medskip

\noindent \textbf{Fact 4.12}\label{R4} \textit{The semiring} $\mathbf{B}_{4}$%
\textit{, defined on the chain }$0<1<$\textit{\ }$2<3$\textit{\ in Example
3.7,} \textit{is not a CI-semiring.\medskip }

\noindent \textit{Proof \ }We shall show that the regular semimodule $_{%
\mathbf{B}_{4}}\mathbf{B}_{4}\in |_{\mathbf{B}_{4}}\mathcal{M}|$ is not
injective. Indeed, consider the semimodule $M\in |_{\mathbf{B}_{4}}\mathcal{M%
}|$ defined on the chain $0<a<b<c<d$ as follows:%
\begin{eqnarray*}
0m &=&0\text{ and }1m=m\text{ for all }m\in M\text{;} \\
\text{ }2a &=&a\text{, }2b=c\text{, }2c=c\text{, }2d=d\text{, }3a=3b=3c=3d=d%
\text{.}
\end{eqnarray*}%
Consider the subsemimodule $K=\{0,b,c,d\}\leq $ $\ _{\mathbf{B}_{4}}M$ and
the homomorphism $\varphi :K\rightarrow \mathbf{B}_{4}$ such that $0\mapsto
0 $, $b\mapsto 1$, $c\mapsto 2$ and $d\mapsto 3$. There is no extention $%
\widetilde{\varphi }:M\longrightarrow \mathbf{B}_{4}$ of the homomorphism $%
\varphi $: Indeed, if such an extention exists, we would have $\widetilde{%
\varphi }(a)\leq \varphi (b)=1$; then, for $\widetilde{\varphi }(a)=0\ $%
would imply $3=\varphi (d)=\varphi (3a)=3\widetilde{\varphi }(a)=0$ and $%
\widetilde{\varphi }(a)=1$ would imply $1=\widetilde{\varphi }(a)=\widetilde{%
\varphi }(2a)=2\widetilde{\varphi }(a)=2$, yielding a contradiction in both
cases. Therefore, the semimodule $_{\mathbf{B}_{4}}\mathbf{B}_{4}$ is not
injective and, consequently, $\mathbf{B}_{4}$ is not a CI-semiring. \textit{%
\ \ \ \ \ \ }$_{\square }$\medskip

\noindent \textbf{Fact 4.13}\label{B(3,1)} \textit{The semiring }$%
B(3,1)=(\{0,1,2\},\oplus ,\odot )$\textit{\ with the operations }$a\oplus b%
\overset{def}{=}\min (2,a+b)$\textit{\ and }$a\odot b\overset{def}{=}\min
(2,ab)$\textit{\ is not a CI-semiring.\medskip }

\noindent \textit{Proof \ }Clearly, $B(3,1)$ is a commutative zerosumfree
anti-bounded semiring (see also \cite[Example 1.8]{Gol1999}). Extend the
monoid $(\{0,1,2\},+,0)$ to a commutative monoid $M=%
\{0,1,2,a_{1},a_{2},a_{3},b_{1},b_{2},b_{3}\}$ with addition
\textquotedblleft $+$\textquotedblright\ defined as follows: $%
a_{i}+a_{i}=a_{i}+b_{i}=b_{i}+b_{i}=b_{i},$ $a_{i}+a_{j}=1$ for $i\neq j$,
and $x+y=2$ for all non-zero $x,y\in M$ in any other case. One can readily
verify that actually $M$ can be naturally considered as a $B(3,1)$%
-semimodule containing $_{B(3,1)}B(3,1)\leq $ $\ _{B(3,1)}M$ $\in |_{B(3,1)}%
\mathcal{M}|$ as a subsemimodule. However, we shall show that $%
_{B(3,1)}B(3,1)$ is not a retract of $_{B(3,1)}M$.

Suppose that there exists a homomorphism $\varphi :M\rightarrow B(3,1)$
which extends $1_{B(3,1)}$. Then, $\varphi (a_{i})=0$ or $\varphi (a_{i})=1$
for each $i\in \{1,2,3\}$: Indeed, otherwise there exists $a_{i}$ such that $%
\varphi (a_{i})=2$, and $1=\varphi (1)=\varphi (a_{i}+a_{j})=\varphi
(a_{i})\oplus \varphi (a_{j})=2\oplus \varphi (a_{j})=2$ for each $j\neq i$.
If $\varphi (a_{j})=0=\varphi (a_{k})$ for at least two different indices $%
j\ $and $k$, then we have a contradiction: $1=\varphi (1)=\varphi
(a_{j}+a_{k})=\varphi (a_{j})\oplus \varphi (a_{k})=0\oplus 0=0$; if $%
\varphi (a_{j})=1=\varphi (a_{k})$ for at least two different indices $j\ $%
and $k$, we have a contradiction: $1=\varphi (1)=\varphi
(a_{j}+a_{k})=\varphi (a_{j})\oplus \varphi (a_{k})=1\oplus 1=2$. Thus,
there is no $\varphi :M\rightarrow B(3,1)$ extending $1_{B(3,1)}$;
therefore, the semimodule $_{B(3,1)}B(3,1)$ is not injective and $B(3,1)$ is
not a CI-semiring. \textit{\ \ \ \ \ \ }$_{\square }$\medskip

\noindent \textbf{Proposition 4.14}\label{+-regular} \textit{Every} \textit{%
anti-bounded left (right) CI-semiring is an additively regular
semiring.\medskip }

\noindent \textit{Proof \ }Let $S$ be an anti-bounded left (right)
CI-semiring. By Proposition 4.1, $S=R\oplus T$\textit{\ }with $R$\ and $T$ a
semisimple ring and a zerosumfree left (right) CI-semiring, respectively. It
is easy to see that $R\ $and $T$,$\ $ being factor semirings of $S$ (see
also \cite[Lemma 5.2]{KN2011}), are a left (right) CI-ring and zerosumfree
anti-bounded CI-semiring, respectively. Therefore, it is enough to show that
the semiring $T$ is an additively regular semiring.

Suppose that a zerosumfree anti-bounded left (right) CI-semiring $T$ with
the multiplicative identity $1\in T$ is not additively regular, \textit{i.e.,%
} $1+x+1\neq 1$ for all $x\in T$. Consider the congruence $\rho \in \mathrm{%
Cong}(T)$ defined as follows: $a\rho b$ iff $a=b$ or there exist $x,y\in T$
such that $a=2+x$ and $b=2+y$ (notice that $2:=1+1\neq 1$ by our
hypothesis). $T/\rho $ is a left (right) CI-semiring. On the other hand, it
is easy to see that $T/\rho \simeq B(3,1)$ which, by Fact 4.13, is not a
CI-semiring. Therefore, $T$ is an additively regular semiring. \textit{\ \ \
\ \ \ }$_{\square }$\medskip

The next result describes all CI-semirings among additively idempotent
anti-bounded semirings.\medskip

\noindent \textbf{Proposition 4.15}\label{ab-B-R3} \textit{An additively
idempotent anti-bounded semiring} $S$ \textit{is a left (right) CI-semiring
iff} $S\simeq \mathbf{B}$ or $S\simeq \mathbf{B}_{3}$.\medskip

\noindent \textit{Proof \ }$\Longrightarrow $. Since $S$ is additively
idempotent, the additive monoid $(S,+,0)$ can be partially ordered by
setting $x\leq y$ iff $x+y=y$. Also, by \cite[Theorem 2.2]{ASW1998}, $S$ is
a regular semiring, and hence, for each nonzero $a\in S$ there exists a
nonzero $x\in S$ such that $axa=a$. Since $S$ is anti-bounded, $x=1+s,$ $%
a=1+s^{\prime }$ for some $s,s^{\prime }\in S$ and, therefore, $%
a=axa=a(1+s)a=a^{2}+asa\geq a^{2}=a(1+s^{\prime })=a+as^{\prime }\geq a$
and, hence, $a^{2}=a$. Similarly, $ab\geq a$ and $ab\geq b$ and, hence, $%
ab\geq a+b$ for all nonzero $a,b\in S$. On the other hand, $%
a+b=(a+b)^{2}=a^{2}+ab+ba+b^{2}=a+ab+ba+b\geq ab$ and, hence, $ab=a+b$. In
particular, $a=a+1$ and, hence, $1\leq a$ for all $0\neq a\in S$. Next, by
Proposition 4.2, $S$ has an infinite element $\infty $, \textit{i.e.}, $%
a+\infty =\infty $ for all $a\in S$, and, therefore, there are the following
cases to consider.

If $\infty =1$ or $S=\{0,1,\infty \}$, then one obviously has $%
S=\{0,1\}\simeq \mathbf{B}$ or $S=\{0,1,\infty \}\simeq \mathbf{B}_{3}$,
respectively.

Now suppose that there exists an element $a\in S$ such that $1<a<\infty $.
If there are no other elements in $S$, then $S\simeq \mathbf{B}_{4}$ and, by
Fact 4.12, $S$ is not a CI-semiring.

Finally, if there exists $b\in S\backslash \{a\}$ such that $1<b<\infty $,
we may assume, without loss of generality, that $b\nleq a$ and consider the
congruence $\rho \in \mathrm{Cong}(S)$ defined as follows: $x\rho y$ iff $x=y
$, or $1<x\leq a\,\,$and $1<y\leq a$, or $x\nleq a\,\,$and $y\nleq a$. Then
one can easily verify that $S/\rho \simeq \mathbf{B}_{4}$ and, hence, by
Fact 4.12, it should be simultaneously not a CI-semiring and a CI-semiring.
So, this case is impossible.

$\Longleftarrow $. This follows from Theorem 4.3 and Fact 4.11. \textit{\ \
\ \ \ \ }$_{\square }\medskip $

Our next observation, in fact, provides a powerful method of constructing a
bunch of interesting anti-bounded semirings arising from arbitrary
rings.\medskip

\noindent \textbf{Example 4.16} Let $R=(R,{+},{\cdot },e,1)$ be an arbitrary
ring with zero $e$ and unit $1$. Let $T:=R\cup \{0\}$ and extend the
operations on $R$ to $T$ by setting $0+t=t=t+0$ and $0\cdot t=0=t\cdot 0$
for all $t\in T$. Clearly, $(T,+,\cdot ,0,1)$ is a zerosumfree semiring.
Now, extend the semiring structure on $T$ to a semiring structure on $%
\mathrm{Ext}(R):=T\cup \{\infty \}=R\cup \{0,\infty \}$, where $\infty
\notin T$, by setting $x+\infty =\infty =\infty +\infty =\infty +x$ and $%
x\cdot \infty =\infty =\infty \cdot \infty =\infty \cdot x$ for all $x\in R$%
, and $0\cdot \infty =0=\infty \cdot 0$. It is easy to see that $(\mathrm{Ext%
}(R),+,\cdot ,0,1)$ is, indeed, an anti-bounded zerosumfree semiring. In a
similar fashion, one can naturally extend the structure of every left $R$%
-module $M$ to a structure of an $\mathrm{Ext}(R)$-semimodule on $\mathrm{Ext%
}(M)$.\medskip

\noindent \textbf{Proposition 4.17 }\textit{The cyclic left (right) }$%
\mathrm{Ext}(R)$\textit{-semimodules are, up to isomorphism,} $\{0\}$\textit{%
, }$\{0,\infty \}$\textit{, and }$\{$\textit{\ }$\mathrm{Ext}(\overline{R})$%
\textit{\ }$|$\textit{\ }$\overline{R}=R/I$\textit{, where }$I$\textit{\ is
a left (right) ideal of }$R$\textit{\ }$\}$\textit{.\medskip }

\noindent \textit{Proof \ }Let $C\in |_{\mathrm{Ext}(R)}\mathcal{M}|$ be a
nonzero cyclic left $\mathrm{Ext}(R)$-semimodule, \textit{i.e.}, $C=\mathrm{%
Ext}(R)c\ $for some $c\in C$. If there exists an element $q\in R$ such that $%
qc=\infty c$, then
\begin{equation*}
q^{\prime }c=(q^{\prime }+e)c=qc+(q^{\prime }-q)c=\infty c+(q^{\prime
}-q)c=(\infty +(q^{\prime }-q))c=\infty c
\end{equation*}

\noindent for any $q^{\prime }\in R$; hence $sc=\infty c$ for all $0\neq
s\in \mathrm{Ext}(R)$ and $C$ is isomorphic to the $\mathrm{Ext}(R)$%
-semimodule $\{0,\infty \}$.

Otherwise, for every $q\in R$, we have $qc\neq \infty c$, and so $qc\neq 0$:
Indeed, $qc=0$ implies $\infty c=\infty qc=0$, and we get a contradiction: $%
0=\infty c=(\infty +1)c=\infty c+c=c$. Thus, $C=\{0\}\cup T\cup \{\infty c\}$
where $T=\{qc\,|\,q\in R\}$, is a cyclic left $R$-module. Whence, $C=\mathrm{%
Ext}(\overline{R})$ where $\overline{R}=R/I$, for some left ideal $%
I\subseteq R$. \textit{\ \ \ \ \ \ }$_{\square }\medskip $

Our next result gives a characterization of semirings $\mathrm{Ext}(R)$, $R$
is a ring, that are CI-semirings.\medskip

\noindent \textbf{Theorem 4.18}\label{Ext(R)} \textit{For a ring }$R=$ $%
(R,+,{\cdot},e,1)$\textit{, the semiring} $\mathrm{Ext}(R)$
\textit{is a left (right) CI-semiring iff }$R$ \textit{is a
semisimple ring.\medskip }

\noindent \textit{Proof \ }$\Longrightarrow $. Let $M$ be a cyclic left $R$%
-module, $A\leq $ $_{R}B$ for $A,B$ $\in |_{R}\mathcal{M}|$, and $\varphi
:A\longrightarrow M$ be an $R$-homomorphism. Then, $\mathrm{Ext}(M)$ is a
cyclic left $\mathrm{Ext}(R)$-semimodule and $\varphi $ induces an $\mathrm{%
Ext}(R)$-homomorphism $\psi :$ $\mathrm{Ext}(A)\longrightarrow \mathrm{Ext}%
(B)$. Since $\mathrm{Ext}(R)$ is a left CI-semiring, the latter can be
extended to $\mathrm{Ext}(R)$-homomorphism $\widetilde{\psi }:\mathrm{Ext}%
(B)\longrightarrow \mathrm{Ext}(M)$. It is easy to see that the restriction $%
\widetilde{\varphi }:=\widetilde{\psi }_{\mid _{B}}:B\longrightarrow M$ is
an $R$-homomorphism extending $\varphi $. Consequently, $R$ is a CI-ring
and, hence, by \cite[p. 649]{Oso1964} (see also \cite[Theorem 1.2.9]{Lam2001}%
), $R$ is a semisimple ring.

$\Longleftarrow $. Let $R=$ $(R,+,\cdot,e,1)$ be a semisimple ring, $%
\mathrm{Ext}(R)=$ $R\cup \{0,\infty \}$, $M=\{\overline{0}\}\cup \overline{R}%
\cup \{\overline{\infty }\}$ a cyclic $\mathrm{Ext}(R)$-semimodule with $%
\overline{R}=R/I$, where $I$\ is a left (right) ideal of $R$, and $_{\mathrm{%
Ext}(R)}M\leq $ $_{\mathrm{Ext}(R)}T\in |_{\mathrm{Ext}(R)}\mathcal{M}|$.
Since $\mathrm{Ext}(R)$ is an additively regular semiring, an additive
reduct of every $\mathrm{Ext}(R)$-semimodule is a commutative inverse monoid
as well. Let $M=[Y;G_{\alpha },\varphi _{\alpha ,\beta }]$ and $%
T=[Z;H_{\alpha },\psi _{\alpha ,\beta }]$ be \textit{Clifford representations%
} of the monoids $(M,+,0)$ and $(T,+,0)$\textit{\ }(see \cite[Theorem 4.11]%
{CP1961} or \cite[Theorem II.2.6]{Pet1984}, also \textit{cf. }\cite[p. 125
and Proposition 2.4]{Kat1997}), respectively. As usual, it is convenient to
identify elements of the semilattices $Y$ and $Z$ in the Clifford
representations with the zeros of the corresponding abelian groups. Then, it
is clear that $Y=\{\overline{0},\overline{e},\overline{\infty }\}\cong
\mathbf{B}_{3}$, $Y\subseteq Z$ and $\psi _{\alpha ,\beta }(a)=a+\beta $ for
any $a\in H_{\alpha }$. Also, it is easy to see that $a+b\in H_{\alpha
+\beta }$ and $sa\in H_{s\alpha }$ for every $a\in H_{\alpha }$, $b\in
H_{\beta }$, and $s\in \mathrm{Ext}(R)$; in particular, since $r\alpha
=\alpha $ for all $r\in R$ and $\alpha \in Z$, one has $\psi _{\alpha ,\beta
}(ra)=ra+\beta =ra+r\beta =r(a+\beta )=r\psi _{\alpha ,\beta }(a)$ for all $%
r\in R$ and $a\in H_{\alpha }$, \textit{i.e.}, $\psi _{\alpha ,\beta }$ are $%
R$-homomorphisms. On the semilattices $Y$ and $Z$, of course, there exists a
natural partial ordering defined as follows: $\alpha \leq \beta $ iff $%
\alpha +\beta =\beta $.

Since $R$ is a semisimple ring, by \cite{Oso1964} (see also \cite[Theorem
1.2.9]{Lam2001} or \cite[Corollary 6.47)]{Lam1999}), there exists an $R$%
-homomorphism $\Theta :H_{\overline{e}}\longrightarrow G_{\overline{e}}$
extending the identity $R$-isomorphism $1_{G_{\overline{e}}}:$ $G_{\overline{%
e}}\longrightarrow G_{\overline{e}}$. We shall show that the map $\Theta
^{\ast }:T\longrightarrow M$ for $a\in H_{\alpha }$, defined by%
\begin{equation*}
\Theta ^{\ast }(a):=\left\{
\begin{tabular}{ll}
$\overline{0}$, & $\text{if }\infty \alpha \text{ }\leq \text{ }\overline{e}$%
, \\
$\Theta \psi _{\alpha ,\beta }$, & if $\alpha \leq \text{ }\overline{e}$ $\&$
$\infty \alpha \text{ }\nleq \text{ }\overline{e}$, \\
$\overline{\infty }$ & in all other cases%
\end{tabular}%
\right. \text{,}
\end{equation*}%
\noindent \noindent is an $\mathrm{Ext}(R)$-homomorphism extending the
identity $\mathrm{Ext}(R)$-homomorphism $1_{M}:$ $M\longrightarrow M$, and
therefore, $M$ is a retract of $T\in |_{\mathrm{Ext}(R)}\mathcal{M}|$.

First, considering all possible cases, we need to verify that $\Theta ^{\ast
}(a+b)=\Theta ^{\ast }(a)+\Theta ^{\ast }(b)$ for all $a\in H_{\alpha }$ and
$b\in H_{\beta }$.

If $\Theta ^{\ast }(a)=\overline{\infty }$, then $\alpha $ $\nleq $ $%
\overline{e}$ and, hence, $a+\beta \nleq $ $\overline{e}$; whence $\Theta
^{\ast }(a+b)=\overline{\infty }=\overline{\infty }+\Theta ^{\ast
}(b)=\Theta ^{\ast }(a)+\Theta ^{\ast }(b)$. The case when $\Theta ^{\ast
}(b)=\overline{\infty }$ can be similarly justified.

If $\Theta ^{\ast }(a)$, $\Theta ^{\ast }(b)\in G_{\overline{e}}$, then $%
\alpha \text{ }\leq \text{ }\overline{e}$, $\beta \text{ }\leq \text{ }%
\overline{e}$ and, hence, $\alpha \text{ }+\beta \text{ }\leq \text{ }%
\overline{e}$ and $\infty (\alpha $ $+\beta )$ $\nleq $ $\overline{e}$;
whence, $\Theta ^{\ast }(a+b)=\Theta (a+b+\overline{e})=$ $\Theta (a+%
\overline{e}+b+\overline{e})=\Theta (a+\overline{e})+\Theta (b+\overline{e}%
)=\Theta ^{\ast }(a)+\Theta ^{\ast }(b)$.

If $\Theta ^{\ast }(a)=\overline{0}=\Theta ^{\ast }(b)$, then $\infty \alpha
\text{ }\leq \text{ }\overline{e}$, $\infty \beta \text{ }\leq \text{ }%
\overline{e}$, and, hence, $\infty (\alpha \text{ }+\beta )\leq \text{ }%
\overline{e}$; whence, $\Theta ^{\ast }(a+b)=\overline{0}=\overline{0}+%
\overline{0}=\Theta ^{\ast }(a)+\Theta ^{\ast }(b)$.

If $\Theta ^{\ast }(a)=\overline{0}$ and $\Theta ^{\ast }(b)\in G_{\overline{%
e}}$, then $\infty \alpha \text{ }\leq \text{ }\overline{e}$, $\beta \leq $ $%
\overline{e}$, $\infty \beta \text{ }\nleq \text{ }\overline{e}$, and,
therefore, $\alpha $ $+\beta \leq \infty \alpha \text{ }+\beta =$ $\beta
\text{ }\leq \text{ }\overline{e}$, and $\infty (\alpha $ $+\beta )$ $\nleq $
$\overline{e}$. Clearly, $\infty a=\infty \alpha $ and so $a+b+\overline{e}=$
$a+b+\overline{e}+\infty \alpha =a+b+\overline{e}+\infty a=b+\overline{e}%
+(1+\infty )a=b+\overline{e}+\infty a=b+\overline{e}$. Therefore, $\Theta
^{\ast }(a+b)=\Theta (a+b+\overline{e})=\Theta (b+\overline{e})=\overline{0}%
+\Theta (b+\overline{e})=\Theta ^{\ast }(a)+\Theta ^{\ast }(b)$. Of course,
the case $\Theta ^{\ast }(b)=\overline{0}$ and $\Theta ^{\ast }(a)\in G_{%
\overline{e}}$ is justified in a similar way.

Now, considering all possible cases, we need to verify that $\Theta ^{\ast
}(sa)=s\Theta ^{\ast }(a)\ $for all $a\in H_{\alpha }$ and $s\in \mathrm{Ext}%
(R)$, where, certainly, we may assume $s\neq 0$.

So, if $\Theta ^{\ast }(a)=\overline{0}$, then $\infty (s\alpha )=\infty
\alpha \leq \text{ }\overline{e}$ and so $\Theta ^{\ast }(sa)=\overline{0}=s%
\overline{0}=s\Theta ^{\ast }(a)$.

If $\Theta ^{\ast }(a)\in G_{\overline{e}}$ and $s\in R$, then $\Theta
^{\ast }(sa)=s\Theta ^{\ast }(a)$ is true since the composite $\Theta \psi
_{\alpha ,\overline{e}}$ is obviously an $R$-homomorphism as well.

If $\Theta ^{\ast }(a)\in G_{\overline{e}}$ and $s=\infty $, then $\infty a%
\text{ }\in $ $H_{\infty \alpha }$ and $\infty \alpha \text{ }\nleq \text{ }%
\overline{e}$ and so $\Theta ^{\ast }(\infty a)=\overline{\infty }=\infty
\Theta ^{\ast }(a)$.

If $\Theta ^{\ast }(a)=\overline{\infty }$, then $\alpha \text{ }\nleq \text{
}\overline{e}$ and, since $\alpha $ $\leq $ $s\alpha $ for $s\neq 0$, one
gets $s\alpha \nleq \text{ }\overline{e}$, and therefore, $\Theta ^{\ast
}(sa)=\overline{\infty }=s\overline{\infty }=s\Theta ^{\ast }(a)$.

Thus, we have shown that $\Theta ^{\ast }:T\longrightarrow M$ is an $\mathrm{%
Ext}(R)$-homomorphism extending the identity $\mathrm{Ext}(R)$-homomorphism $%
1_{M}:$ $M\longrightarrow M$, and therefore, $M$ is a retract of $T\in |_{%
\mathrm{Ext}(R)}\mathcal{M}|$. Now, taking into consideration that $\mathrm{%
Ext}(R)$ is an additively regular semiring and applying \cite[Theorem 4.2 or
Corollary 4.3]{Kat1997}, we can choose the semimodule $T\in |_{\mathrm{Ext}%
(R)}\mathcal{M}|$ to be injective and, hence, conclude that $M$ is an
injective $\mathrm{Ext}(R)$-semimodule, too. To finish the proof, one needs
only to use Proposition 4.17 and note that the $\mathrm{Ext}(R)$-semimodule $%
\{0,\infty \}$ is obviously a retract of the regular semimodule $_{\mathrm{%
Ext}(R)}\mathrm{Ext}(R)$. \textit{\ \ \ \ \ \ }$_{\square }\medskip $

In the next observation, and as a consequence of the previous theorem, we
obtain a complete description of zerosumfree additively regular anti-bounded
CI-semirings.\medskip

\noindent \textbf{Proposition 4.19}\label{iso} \textit{A zerosumfree
additively regular anti-bounded} \textit{semiring} $S$ \textit{is a left
(right) CI-semiring iff }$S\simeq \mathbf{B}$\textit{, or} $S\simeq \mathbf{B%
}_{3}$\textit{, or} $S\simeq \mathrm{Ext}(R)$ \textit{for some nonzero
semisimple ring} $R$\textit{.\medskip }

\noindent \textit{Proof \ }$\Longrightarrow $. By Proposition 4.1, $S$
contains an infinite element $\infty $. Consider the congruence $\diamond
\in \mathrm{Cong}(S)$ defined as follows: $s\diamond s^{\prime }$ iff $%
ns=s^{\prime }+t$ and $ms^{\prime }=s+t^{\prime }$ for some $n,m\in \mathbb{N%
}$ and $t,t^{\prime }\in S$. By \cite[Lemma 2.2]{Ili2012}, the quotient
semiring $\overline{S}=S/\diamond $ is an additively idempotent, left
(right) CI-semiring. Thus, applying Proposition 4.15, $\overline{S}\simeq
\mathbf{B}$ or $\overline{S}\simeq \mathbf{B}_{3}$, and we will consider
each of these cases.

If $\overline{S}\cong \mathbf{B}$, then$\ $we have $1\diamond \infty \ $and,
hence, $n1=\infty +s=\infty $ for some $n\in \mathbb{N}$ and $s\in S$. Since
$S $ is an additively regular semiring, there exists $x\in S$ such that $%
1=1+x+1 $ and, therefore, $1=n(1+x)+1=n1+nx+1=\infty +nx+1=\infty $. Whence,
$S$ is even an additively idempotent semiring and, clearly, $S\cong
\overline{S}\simeq \mathbf{B}$.

Now, let $\overline{S}\cong \mathbf{B}_{3}$, and $x\in S$ be an additive
inverse of $1\in S$, \textit{i.e.}, $1+x+1=1$. Then, for the element $e:=1+x$
$\in I^{+}(S)$ and all $s\in S$ and $a\in I^{+}(S)$, we have
\begin{equation*}
e^{2}=(1+x)^{2}=1+x+x+x^{2}=1+x(1+1+x)=1+x\cdot 1=1+x=e,
\end{equation*}%
\begin{eqnarray*}
se &=&s+sx=(1+1+x)s+sx=s+s+xs+sx=xs+s(1+1+x)= \\
xs+s &=&es\ \text{and} \\
a &=&a\cdot 1=a(1+1+x)=a+a+ax=a+ax=a(1+x)=ae\text{.}
\end{eqnarray*}%
It is clear that $e\neq 0$, $e\neq \infty \ $and the restriction $\pi
_{|_{I^{+}(S)}}$ of the natural surjection $\pi :S\longrightarrow \overline{S%
}$ is an injection. Therefore, $I^{+}(S)=\{0,e,\infty \}$.

As was shown above for $a=1$, it can be shown that $a\diamond \infty
$ implies $a=\infty $, and therefore, the equivalence classes
$0^\diamond$ and $\infty^\diamond$ are $\{0\}$ and $\{\infty \}$,
respectively. Thus,
$S=\{0\}\cup R\cup \{\infty \}$, where $R\ $ is the equivalence class $%
e^{\diamond }$. If $|R|=1$, then it is easy to see that $S\cong \mathbf{B}%
_{3}$. Therefore, we have only to consider the case when $|R|>1$. Obviously,
$R$ is closed under addition and multiplication. Hence, for any $a\in R$ we
have $ae\in R\cap I^{+}(S)$ and so $ae=e$ and, consequently, $%
a+e=a+ae=a(1+e)=a\cdot 1=a$ and $a+ax=a(1+x)=ae=e$. Whence, $(R,{+},{\cdot }%
,e,1)$ is a ring and $S=\mathrm{Ext}(R)$ and, by Theorem 4.18, $R$ is a
semisimple ring.

$\Longleftarrow $. This follows from Proposition 4.15 and Theorem 4.18.
\textit{\ \ \ \ \ \ }$_{\square }\medskip $

Applying Propositions 4.1, 4.15 and 4.19, we obtain a complete
characterization of anti-bounded CI-semirings generalizing Osofsky's
celebrated characterization of semisimple rings (\cite{Oso1964}, see also
\cite[Theorem 1.2.9]{Lam2001} and \cite[Corollary 6.47)]{Lam1999})\medskip :

\noindent \textbf{Theorem 4.20}\label{ab-CI} \textit{An anti-bounded}
\textit{semiring }$S$ \textit{is a left (right) CI-semiring iff }$S$ \textit{%
is one of the following semirings:}

\textit{(1) }$S$ \textit{is a semisimple ring;}

\textit{(2)} $S\simeq \mathbf{B},$ \textit{or} $S\simeq \mathbf{B}_{3},$
\textit{or} $S\simeq \mathrm{Ext}(R)$ \textit{for some nonzero semisimple
ring} $R$\textit{;}

\textit{(3)} $S=R\oplus T$\textit{,} \textit{where} $R$ \textit{is a
semisimple ring and} $T$ \textit{is isomorphic to} $\mathbf{B},$ \textit{or}
$\mathbf{B}_{3},$ \textit{or} $\mathrm{Ext}(R^{\prime })$ \textit{for some
nonzero semisimple ring} $R^{\prime }$\textit{.\medskip }

Finally, as an application of Theorem 4.20, we solve \cite[Problem 3]%
{KNT2011} and \cite[Problem 5]{KNZ} in the class of anti-bounded
CI-semirings.\medskip

\noindent \textbf{Corollary 4.21}\label{perfect} \textit{For an anti-bounded
semiring} $S,$ \textit{the following conditions are equivalent:}

\textit{(i)} $S$ \textit{is a congruence-simple left (right) CI-semiring;}

\textit{(ii)} $S$ \textit{is an ideal-simple left (right) CI-semiring;}

\textit{(iii)} $S\cong M_{n}(D)$ \textit{for some division ring} $D$ \textit{%
and} $n\geq 1,$ \textit{or} $S\cong \mathbf{B}$.

\end{document}